\documentclass{amsart}
\usepackage[utf8]{inputenc}

\usepackage[english]{babel}

\usepackage{amssymb}
\usepackage{amsmath}
\usepackage{amsthm}
\usepackage{enumerate}
\usepackage{cite}
\usepackage{xcolor}
\usepackage[final]{listings}
\usepackage[section]{algorithm}

\makeindex

\title[Beyond Kruskal's bound]
{Identifiability beyond Kruskal's bound
for symmetric tensors of degree $4$}
\date{}

\newcommand{\C}{\mathbb{C}}
\newcommand{\Z}{\mathbb{Z}}
\newcommand{\Pj}{\mathbb{P}}

\newcommand{\N}{\mathbb{N}}

\newcommand{\Oc}{\mathcal{O}}

\newcommand{\vect}[1]{\mathbf{#1}}

\newcommand{\Sec}[2]{\sigma_{#1}({#2})}
\newcommand{\Tang}[2]{\mathrm{T}_{#1} {#2}}

\newcommand{\rank}{\operatorname{rank}}
\newcommand{\imm}{\operatorname{im}}

\newcommand{\refsec}[1]{{section \ref{#1}}}

\newcommand{\refprop}[1]{{Proposition \ref{#1}}}
\newcommand{\reflem}[1]{{Lemma \ref{#1}}}

\newtheorem{thm0}{Theorem}[section]
\newtheorem{prop0}[thm0]{Proposition}
\newtheorem{lemma0}[thm0]{Lemma}
\newtheorem{coro0}[thm0]{Corollary}
\theoremstyle{definition}
\newtheorem{defn0}[thm0]{Definition}

\newtheorem{exa0}[thm0]{Example}

\newtheorem{rem0}[thm0]{Remark}

\subjclass[2000]{14J70, 14C20, 14N05, 15A69, 15A72}

\author[E.~Angelini]{Elena Angelini}
\address{Dipartimento di Ingegneria dell'Informazione e Scienze Matematiche, Universit\`a di Siena, Italy}
\email{elena.angelini@unisi.it}

\author[L.~Chiantini]{Luca Chiantini}
\address{Dipartimento di Ingegneria dell'Informazione e Scienze Matematiche, Universit\`a di Siena, Italy}
\email{luca.chiantini@unisi.it}
\thanks{The first and second author are members of the Italian GNSAGA-INDAM and are supported by the Italian PRIN 2015 - 
Geometry of Algebraic Varieties (B16J15002000005)}

\author[N.~Vannieuwenhoven]{Nick Vannieuwenhoven}
\address{Department of Computer Science, KU Leuven, Belgium}
\email{nick.vannieuwenhoven@cs.kuleuven.be}
\thanks{The third author was supported by a Postdoctoral\ Fellowship of the 
Research Foundation--Flanders (FWO)}

\begin{document}

\begin{abstract}
We show how methods of algebraic geometry can produce criteria for the identifiability of
specific tensors that reach beyond the range of applicability of the celebrated Kruskal criterion.
More specifically, we deal with the symmetric identifiability of symmetric tensors in $\mathrm{Sym}^4(\C^{n+1})$,
i.e., quartic hypersurfaces in a projective space $\Pj^n$, that have a decomposition
in $2n+1$ summands of rank 1. This is the first case where the reshaped Kruskal criterion 
no longer applies.
We present an effective algorithm, based on efficient linear algebra computations, that checks
if the given decomposition is minimal and unique. The criterion is based on the application
of advanced geometric tools, like Castelnuovo's lemma for the existence of rational
normal curves passing through a finite set of points, and the Cayley--Bacharach condition
on the postulation of finite sets. In order to apply these tools to our situation, 
we prove a reformulation of these results, hereby extending classical results such as
Castelnuovo's lemma and the analysis of Geramita, Kreuzer, and Robbiano, Cayley--Bacharach schemes and their canonical modules, \textit{Trans. Amer. Math. Soc.} 339:443--452, 1993.
\end{abstract}

\maketitle

\section{Introduction}
The aim of this paper is the continuation of the study, started in \cite{COttVan17b}, of conditions which imply the
identifiability of symmetric tensors, in a numerical range where the celebrated Kruskal criterion does not apply. Recall that a symmetric tensor $T$ is \emph{identifiable} if there exists a unique decomposition $T=T_1+\dots+T_r$ with a minimal number of symmetric rank-$1$ terms, up to scaling and reordering of the summands. This decomposition is called a \textit{symmetric tensor rank decomposition} or \textit{Waring decomposition}.

Beyond its theoretical interest, identifiability plays a central role in many applications of symmetric tensors. An important class of applications is found in algebraic statistics and machine learning.
Indeed, the parameters of several latent variable models, including topic models, latent Dirichlet allocation, and hidden Markov models, can be recovered from the unique Waring decomposition of a symmetric tensor $T$ that is associated with the model, as shown in \cite{AnandkumarGeHsuKakadeTelgarsky14,AllmanMatiasRhodes09,SKLV2017}. Since one wishes to recover and interpret the parameters of the model uniquely, it is important to verify that the symmetric decomposition (computed by some numerical algorithm) is unique. 
In 1977, Kruskal \cite{Kruskal77} determined what is now still the most popular criterion for testing the identifiability of a
specific decomposition of a given, general tensor. 
Kruskal showed that a given decomposition is identifiable if its length is smaller than some numerical condition on the \emph{Kruskal ranks} associated with the decomposition. Applying Kruskal's result to a reshaped tensor is considered to be a state-of-the-art effective criterion for identifiability by \cite{COttVan17b}.

Despite its popularity, there are three limitations to Kruskal's test in the symmetric setting. First, Kruskal's criterion verifies whether the tensor has no other tensor rank decompositions. It is possible, in principle, that a symmetric tensor has only one symmetric decomposition but several tensor rank decompositions of the same length. Indeed, Shitov \cite{Shitov2017} recently provided a counterexample to Comon's conjecture, which entails that this may happen. Second, the length $r$ of the decompositions of which Kruskal's criterion can prove identifiability is much smaller than the range wherein generic symmetric identifiability holds \cite{COttVan17a}. Third, Derksen proved in \cite{Derksen13} that the numerical condition on the Kruskal ranks is sharp in the sense that if $s$ is the largest rank allowed by the numerical condition in Kruskal's criterion, then there exist unidentifiable tensors of rank $s+1$ that still satisfy the numerical condition. 

\begin{rem0}
The construction in \cite{Derksen13} also applies to the symmetric case. Indeed, taking all Kruskal ranks equal to the maximum value, the example constructed in \cite{Derksen13} becomes a symmetric tensor decomposition $T = T_1 + \cdots + T_{s+1}$ that admits another, distinct symmetric tensor decomposition $T = T_1' + \cdots + T_{s+1}'$; hence, \cite{Derksen13} also proved that Kruskal's criterion is sharp in the symmetric setting.
\end{rem0}

The foregoing reasons motivated a further study of the specific identifiability of symmetric tensors, as
in \cite{BallC12, COttVan17b,MassaMellaStagliano}. However, results beyond the (reshaped) Kruskal criterion are sparse. Moreover, 
some of these criteria require the use of general computer algebra algorithms which rapidly become ineffective for high-dimensional varieties.
In Proposition 6.3 of \cite{COttVan17b}, a criterion for the identifiability of quartics in $\Pj^3$ (symmetric tensors of type $4\times 4 \times 4 \times 4$), having a decomposition with $7$ summands was given; this is exactly the first value beyond the range of the reshaped Kruskal's criterion.
The goal of this paper is to extend the analysis to symmetric tensors in $S^{4} \C^{n+1} $, for any $n$.
For such tensors the reshaped Kruskal criterion can prove the identifiability
only for decompositions with $2n$ or less summands. We propose here a test that
verifies the identifiability even for decompositions of length $2n+1$. Moreover, the test is
\textit{effective}, in the sense of \cite{COttVan17b}: it will give a positive answer except on a set of measure $0$ in the variety of tensors decomposed with $2n+1$ or less summands.

The idea of the test is based on the following observation. If $T\in S^{4} \C^{n+1}$ has a decomposition
with $2n+1$ summands it may happen that the decomposition is not unique, even if its associated Kruskal ranks have the maximum value $n+1$, which is consistent with \cite{Derksen13}. However, in this case, we prove the existence of a positive-dimensional family of decompositions for $T$, which includes the given one.
Specifically, we prove that when $T$ has two different decompositions of length $2n+1$, then the pre-images of the two decompositions in the Veronese map determine a finite subset of points $Z$ in a projective space $\Pj^n=\Pj(\C^{n+1})$ which lies in a rational normal curve. This rational normal curve in $\Pj^n$ induces the existence of the positive-dimensional family of decompositions for the tensor $T$, which can be detected using only linear algebra. These insights yield the new criterion (see \refsec{sec_the_test}).

We notice that the failure of identifiability for tensors of rank $2n+1$ in $S^{4} \C^{n+1}$ is caused by the existence of a low dimensional variety---the rational normal curve---that contains the pre-images of the points of the decomposition and forces the existence of other decompositions of the same tensor. The link between the failure of identifiability and the existence of subvarieties of positive dimension containing the decompositions is somehow familiar in the study of tensors. This phenomenon can occur for generic tensors, where the varieties that force the existence of many decompositions of $T$ are the \textit{contact varieties}; see \cite{CCi06}. In the case of symmetric tensors of very low rank, it was shown in \cite{BallC13} that unidentifiable cases occur only if the decompositions are contained in a positive-dimensional subvariety.
In this paper, we prove that the same fact holds for specific tensors of rank $2n+1$ in  $S^{4} \C^{n+1}$, and we expect that a similar characterization of unidentifiable specific tensors can be proved in other cases which lie just beyond Kruskal's range. 

The basic tool in our analysis is provided by the study of the geometry of finite sets $Z$ in
the projective space $\Pj^n$. We will perform the analysis by means of classical methods in algebraic
geometry, essentially related to the Hilbert function of $Z$. However, we cannot plainly use the large body of classical and modern results on finite sets in projective spaces. The reason is that when we have
a decomposition of $T$, and we want to exclude the existence of a second decomposition, then we 
argue on the pre-images of the two decompositions, which determine two sets $A,B\subset\Pj^n$.
In order to achieve our result, i.e., that $Z=A\cup B$ lies in a rational normal curve, we can control
only the geometry of $A$, as we know nothing about the hypothetical set $B$; in particular, we cannot place assumptions on the geometry of $B$. For this reason, we need to produce refinements of well-known geometric results, such as Castelnuovo's Lemma (see \reflem{lemma:Cclass}), in which we sharpen the hypothesis on the generality of the position of the points in $Z$ in \reflem{lemma:Cnew}. Similarly, we introduce a \textit{relative} version of the Cayley-Bacharach condition (see Definition \ref{def:CB}), and extend a result of Geramita, Kreuzer and Robbiano (see Theorem \ref{GKRext} below).

We hope that our analysis can be of independent interest in the theory of finite sets in projective spaces. We also believe that it can support the idea that geometric results
on the geometry of sets of points can produce interesting consequences for the theory of symmetric tensors. We also strongly believe that further analyses of the same type can provide new applications of algebraic geometry methods in the study of tensors, as well as stimulate the research on the Hilbert functions of finite sets, by indicating which refinements of known results could produce non-trivial applications to tensor analysis.

The rest of this article is structured as follows. In the next section some elementary results about the Hilbert function of finite sets are recalled. Kruskal's identifiability criterion of tensors is recalled in \refsec{sec:Kruskal}. We then investigate, in \refsec{sec:CB&res}, the Hilbert function of sets with the Cayley--Bacharach property. In \refsec{sec:Crev}, the assumptions in the classic Castelnuovo Lemma are relaxed. Finally, we apply the results from aforementioned sections to the identifiability of fourth-order symmetric tensors whose rank is one higher than the range in which the (reshaped) Kruskal criterion applies.

\paragraph{\textbf{Acknowledgements.}} The authors would like to thank Juan Migliore, Giorgio Ottaviani, and Maria Evelina Rossi for several fruitful discussions on the topics of the present research.

\section{Preliminaries}\label{sec:notation}

\subsection{Notation}
Let $ T $ be a homogeneous polynomial in $ n+1 $ variables of degree $ d $ over $ \C $, i.e., $ T \in S^{d} \C^{n+1} $. 
$ T $ is associated to an element of $ \Pj(S^{d} \C^{n+1}) $, which by abuse of notation we denote by $T$.

For any $ m \in \N $, let $\Pj^m=\Pj^m _\C$ be the $ m $-dimensional complex projective space and let 
$\nu_{d}:\Pj^n\to \Pj^N$ be the Veronese embedding of $ \Pj^n $ of degree $ d $, where $N = \binom{n+d}d - 1 $.
Let $ A \subset \Pj^n $ be a finite set. We denote by $ \ell(A) $ the cardinality of $ A $ and we define 
$$ \nu_{d}(A) = \{\nu_{d}(P_{1}), \ldots, \nu_{d}(P_{\ell(A)}) \} \subset \Pj^{N}, $$
where $ P_{i} \in A $.  
\smallskip

With the above notations we give the following definitions.

\begin{defn0}
Let $ A \subset \Pj^n $ be a finite set. $ A $ \emph{computes} $ T $ if $ T\in \langle \nu_{d}(A) \rangle$, 
the linear space spanned by the points of $ \nu_{d}(A) $.
\end{defn0}

Recall that the Kruskal rank of a finite set $A \subset \Pj^n$ is defined as the maximum value $k$ such that all subsets of $k$ points from $A$ are linearly independent. By definition, the maximum value for the Kruskal rank of $A$ is thus $\min\{ \ell(A), n+1\}$, which is also the generic value.

\begin{defn0}
A finite set $ A \subset \Pj^n $ is in \emph{linear general position} (LGP) if the Kruskal rank of $A$ is maximal, i.e., equal to $\min\{\ell(A), n+1\}$.
This implies that for any $a\leq n+1$, any subset of $A$ of cardinality $a$ is linearly independent. 
\end{defn0}

\begin{defn0}
Let $ A \subset \Pj^n $ be a finite set which computes $ T $. $ A $ is \emph{minimal} if we cannot find a proper subset $ A' $ of $ A $ 
such that $T \in \langle\nu_{d}(A')\rangle $.
\end{defn0}

\begin{rem0}\label{rem:indep}
If $ A \subset \Pj^n $ is a finite set that computes $T$ and satisfies the minimality property, then the points of $ \nu_{d}(A) $ are 
linearly independent, i.e.,
$$ \dim(\langle\nu_{d}(A)\rangle) = \ell(A) -1.  $$
\end{rem0}

Recall that a \emph{rational normal curve} is a curve $\Gamma \subset \Pj^n $ corresponding to the Veronese embedding of $ \Pj^{1} $ of degree $ n $.
Rational normal curves are the only irreducible curves of degree $n$ in $\Pj^n$.

\begin{rem0}\label{rem:rncLGP}
If $ Z $ is a finite subset of a rational normal curve $ \Gamma \subset \Pj^{n} $, then $ Z $ is in LGP and satisfies $ Dh_{Z}(j) \leq n $ for $j>0$; see \cite{Harris} at the bottom of page 10.
\end{rem0}

\begin{rem0}\label{rem:rncsec}
It is classically known that curves are never defective, i.e., their secant varieties always have the expected dimension.
Thus, if $\Gamma \subset \Pj^n $ is a rational normal curve, then the dimension of the $ k $-secant variety 
$ \sigma_{k}(\Gamma) $ is the expected
$$ \min\{n,2k-1\}. $$
\end{rem0}

\subsection{The Hilbert function of finite sets and its difference}\label{sec:hilb}

\begin{defn0}
Let $ Y \subset \C^{n+1} $ be a finite set of cardinality $ \ell $ and let $ j \in \N $. The \emph{evaluation map 
of degree $ j $ on $ Y $} is the linear map
$$ ev_{Y}(j): S^{j} \C^{n+1} \to \C^\ell $$ 
which sends $ F \in S^{j} \C^{n+1} $ to the evaluation of $ F$ at the points of $ Y $. 
\end{defn0}

Let $ Z \subset \Pj^n $ be a finite set. We give the following definition.

\begin{defn0}
The \emph{Hilbert function} of $ Z $ is the map
$$ h_{Z} : \Z \to \N $$
such that $ h_{Z}(j) = 0 $ for $ j < 0 $ and $ h_{Z}(j) = \rank(ev_{Y}(j)) $ for $ j \geq 0 $, where $ Y\subset \C^{n+1}$ 
is any set of $ (n+1) $-tuples representing the elements of $ Z $, i.e., $ [Y] = Z $. 
\end{defn0}

\begin{rem0}
One verifies that while $ ev_{Y}(j) $ depends on the choice of $Y$, such that $[Y] = Z$, 
its rank does not depend on this choice. Therefore, $ h_{Z} $ is well defined.  
\end{rem0}

\begin{rem0}\label{rem:conds}
For any $ j \geq 0 $, the value $ h_{Z}(j) $ provides the number of conditions that $ Z $ imposes to 
the elements of $ S^{j} \C^{n+1} $, i.e., $ h_{Z}(j) = \dim (\langle \nu_{j}(Z) \rangle) +1 $. In particular, 
if $ h_{Z}(j) = \ell(Z) $, then we say that \emph{$Z$ imposes independent conditions to forms of degree $j$}.
\end{rem0}

\begin{defn0}
The \emph{first difference of the Hilbert function} $ Dh_{Z} $ of $ Z $ is given by
$$ Dh_{Z}(j) = h_{Z}(j)-h_{Z}(j-1),$$
where $ j \in \Z $.
The set of non-zero values of $Dh_Z$ is called the \emph{h-vector} of $ Z $.
\end{defn0}

We recall some elementary and well-known properties of $ h_{Z} $ and $ Dh_{Z} $ that will be useful throughout the paper.

\begin{lemma0}\label{rem:triv} We have
\begin{enumerate}[(i)]
\item $ Dh_{Z}(j) = 0$ for $ j < 0$;
\item $ h_{Z}(0) = Dh_{Z}(0) = 1$;
\item $ Dh_{Z}(j) \geq 0 $ for all $ j $;
\item $ h_{Z}(j) = \ell(Z) $ for all $ j \gg 0$; 
\item $ Dh_{Z}(j) = 0$ for $ j \gg 0 $ and $ \sum_{j} Dh_{Z}(j) = \ell(Z) $; 
\item $ h_{Z}(i) = \sum_{0\leq j \leq i} Dh_{Z}(j) $; 
\item if $ h_{Z}(j) = \ell(Z) $, then $ Dh_{Z}(j+1) = 0 $.
\item If $ Z' \subset Z $, then, for every $ j \in \Z $, we have $ h_{Z'}(j) \leq h_{Z}(j) $ and
$Dh_{Z'}(j) \leq Dh_{Z}(j).$
\end{enumerate}
\end{lemma0}

The next property is a consequence of the Macaulay maximal growth principle;
see, e.g., section 3 of \cite{BigaGerMig94} for a proof.

\begin{prop0}\label{nonincr}
If for some $j>0$, $Dh_{Z}(j) \leq j$, then
$$ Dh_{Z}(j) \geq Dh_{Z}(j+1). $$
In particular, if  for some $j>0$, $ Dh_{Z}(j)=0 $, then $Dh_{Z}(i)=0$ for all $i\geq j$.
\end{prop0}

\begin{rem0}\label{indepcond}
Notice that if for some $j$ we have  $ Dh_{Z}(j) = 0 $,  then $ h_{Z}(j-1) = h_{Z}(j) $. By Proposition \ref{nonincr},
 for any $ i \geq j $ also $ Dh_{Z}(i) = 0 $, i.e.,
$ h_{Z}(j-1)=h_Z(i)$ for any $ i \geq j $. Therefore, by parts (v) and (vi) of \reflem{rem:triv}, 
$$ \ell(Z) = \sum_{k} Dh_{Z}(k) = \sum_{k=0}^{i} Dh_{Z}(k) = h_{Z}(i).  $$
Thus, $h_Z(j-1)$ is equal to the cardinality of $Z$, i.e., the evaluation map in degree $j-1$ surjects.
In this case, for every $P\in Z$ we can find a form of degree $i$ that vanishes at
$Z\setminus \{P\}$ and does not vanish at $P$.
Therefore, when $h_Z(i)=\ell(Z)$, we will also say that \textit{hypersurfaces of degree $i$ separate the points of $Z$}.
\end{rem0}

\begin{rem0}\label{resto1}
Assume $h_Z(i)=\ell(Z)-1$. Then $h_Z(i+1)>h_Z(i)$, for otherwise, by \refprop{nonincr}, $Dh_Z(j)=0$ for all $j>i$,
thus $h_Z(j)=h_Z(i)$ for all $j>i$, contradicting property (iv) of \reflem{rem:triv}.
Thus, if $h_Z(i)=\ell(Z)-1$ then necessarily $h_Z(i+1)=\ell(Z)$.
\end{rem0}

\begin{rem0}\label{exseq}
There is an alternative way to look at the Hilbert function of a finite set $ Z \subset \Pj^{n} $. Indeed, let $ J_{Z} $ be the 
ideal sheaf of $ Z $ and let $ j \in \N $. Then, we have the short exact sequence
\begin{equation}\label{eq:exs}
0 \longrightarrow J_{Z}(j) \longrightarrow \Oc_{\Pj^{n}}(j) \longrightarrow \Oc_{Z}(j) \longrightarrow 0. 
\end{equation}
Passing to cohomology, \eqref{eq:exs} provides the exact sequence
$$ 0 \longrightarrow H^{0}(J_{Z}(j)) \longrightarrow H^{0}(\Oc_{\Pj^{n}}(j)) \buildrel \rm {\it ev_{Y}(j)} \over 
\longrightarrow H^{0}(\Oc_{Z}(j)) \longrightarrow H^{1}(J_{Z}(j)) \longrightarrow 0, $$
regardless of the choice of $ Y $ with $ [Y] = Z $. Therefore,
\begin{equation}\label{eq:h^0}
h_{Z}(j) =  \binom{j+n}j - \dim(H^{0}(J_{Z}(j))),
\end{equation}
and, by \reflem{rem:triv} (v) and (vi),
\begin{equation}\label{eq:h^1}
 \dim(H^{1}(J_{Z}(j))) = \ell(Z)-h_{Z}(j) = \sum _{i> j} Dh_{Z}(i).
\end{equation}
\end{rem0}

\section{Kruskal's criterion for symmetric tensors}\label{sec:Kruskal}

In this section, we recall the specialization of the celebrated Kruskal criterion for the identifiability of a decomposition to symmetric tensors. In fact, in the context of this paper, we are interested only in the generic case where the Kruskal ranks are maximal. In this case, the highest rank $r$ of which Kruskal's criterion can prove identifiability is maximal. The particular version of Kruskal's Lemma that is of relevance is recalled as the following corollary of the results in \cite{Kruskal77}.

\begin{coro0}[Kruskal \cite{Kruskal77}]\label{thm:Krusk}
Let $ T $ be an $ n_{1} \times n_{2} \times n_{3} $ tensor over $ \C $ with $ n_{1} \geq n_{2} \geq n_{3} \geq 2 $. Assume that $ T = T_{1}+ \cdots + T_{r} $, 
where the $ T_{i} $'s are tensors of rank $ 1 $. Write $T_i=v_{1i}\otimes v_{2i}\otimes v_{3i}$. If the sets
$A_j=\{v_{j1},\dots,v_{jr}\}$ are in LGP and 
$$ r \leq \frac{1}{2} (\min(n_{1},r)+\min(n_{2},r)+\min(n_{3},r)) - 1$$
then $ T $ has complex rank $ r $ and it is identifiable, in the sense that the set $\{T_1,\ldots,T_r\}$ is unique, including multiplicities.
\end{coro0}

A value of $r$ not satisfying the above inequality is said to be \emph{beyond Kruskal's range} of identifiability.

According to Corollary 20 of \cite{COttVan17b}, Kruskal's criterion can be specialized to the case of symmetric tensors, as follows. As before, we only present the following corollary in the generic case where the Kruskal ranks are maximal.

\begin{coro0}[Reshaped Kruskal's criterion for symmetric tensors \cite{COttVan17b}]\label{thm:Kruskalresh} 
Let $ T \in S^{d} \C^{n+1} $ with $ d \geq 3 $ and $n\ge1$. Let $ A \subset \Pj^{n} $ be a finite set of cardinality $r = \ell(A)$ computing $ T $, and let 
$ d_{1}+d_{2}+d_{3} = d $ be a partition of $ d $ such that $ d_{1} \geq d_{2} \geq d_{3} $. If $\nu_{d_{i}}(A)$ is in LGP for $ i=1,2,3 $ and
$$ 
r \leq \frac{1}{2} \bigl( \min\{ \tbinom{d_{1}+n}{d_{1}}, r \} + \min\{ \tbinom{d_{2}+n}{d_{2}}, r \} + \min \{ \tbinom{d_{3}+n}{d_{3}}, r \} \bigr) -1,
$$
then $ T $ has complex rank $ r $ and it is identifiable.
\end{coro0}

\begin{rem0}\label{rem:last}
Direct computations show that for the case $ d = 4 $, which is the core of this paper, the maximum range of applicability is attained 
for $ d_{1} = 2 $ and $ d_{2} = d_{3} = 1 $, so that $ \ell(A) \leq 2n $. Values of $r = \ell(A) > 2n$ are beyond the reshaped Kruskal's range.
\end{rem0}

\section{Identifiability and the Cayley-Bacharach property}\label{sec:CB&res}

For a finite set, we define the Cayley-Bacharach property $CB(i)$ as follows.

\begin{defn0}\label{def:CB}
A finite set $Z\subset \Pj^n$ satisfies the \emph{Cayley-Bacharach property in degree $i$}, abbreviated as $\mathit{CB}(i)$, if for all $P \in Z$, it holds that every form of degree $i$ vanishing at $ Z\setminus\{ P\}$
also vanishes at $P$.
\end{defn0}

\begin{rem0}\label{rem:CBprop}
If $ Z $ satisfies $\mathit{CB}(i)$, then it satisfies $\mathit{CB}(i-1)$ too. Otherwise, one could find 
$ P \in Z $ and a  hypersurface $ F \subset \Pj^{n} $ of degree $ (i-1) $ such that $ Z \setminus \{P\} 
\subset F $ and $ P \notin F $. Therefore, if $ H_{P} \subset \Pj^{n} $ is a hyperplane 
not containing $ P $, then $ F \cup H_{P} \in H^{0}(J_{Z\setminus\{ P\}}(i)) \setminus H^{0}(J_{Z}(i)) $, 
which contradicts the hypothesis.
\end{rem0}

\begin{rem0}
If $ Z $ satisfies $\mathit{CB}(i)$, then for any $ P \in Z $, we have
$$ H^{0}(J_{Z\setminus\{ P\}}(i)) = H^{0}(J_{Z}(i)).$$
It follows from equation \eqref{eq:h^0} and Remark \ref{rem:CBprop} that
\begin{equation}\label{eq:h0}
h_Z(j)=h_{Z\setminus\{P\}}(j) \quad \mbox{ and }\quad Dh_Z(j)=Dh_{Z\setminus\{P\}}(j) \quad \forall j\leq i.
\end{equation}
\end{rem0}

The Cayley-Bacharach property could thus be interpreted as the converse of the separation property introduced in Remark \ref{indepcond}.

\begin{rem0}
From Remark \ref{indepcond} it is clear that if $Z$ satisfies $\mathit CB(i)$,  then  hypersurfaces of degree $i$ cannot separate the points of $Z$. We must namely have
\begin{equation}\label{eq:hl}
h_{Z}(i) < \ell(Z)  \qquad\mbox{or, equivalently}\qquad Dh_Z(i+1)>0.
\end{equation} 
Indeed, if $ h_{Z}(i) = \ell(Z) $, then from \eqref{eq:h0} and the exact sequence
$$ 0 \to H^{0}(J_{Z\setminus\{ P\}}(i)) \to H^{0}(\Oc_{\Pj^n}(i)) \to H^{0}(\Oc_{Z\setminus\{ P\}}(i)) \to H^{1}(J_{Z\setminus\{ P\}}(i)) \to 0,  $$
it follows that 
$$ h^{1}(J_{Z\setminus\{ P\}}(i)) = \ell(Z\setminus\{ P\}) - h_{Z\setminus\{ P\}}(i) = \ell(Z) - 1 - h_{Z}(i) = -1, $$
which is not possible.

We notice that the converse statement is false. For instance, the set $ Z $ consisting of four points in $ \Pj^2 $, three of them aligned,
does not satisfy $\mathit{CB}(1)$, while $h_{Z}(1)<4$.
\end{rem0}

For brevity, we define the following value.
\begin{defn0}
The \emph{socle degree} $ i_{Z} $ of $ Z $ is the maximum $ i $  such that \eqref{eq:hl} holds.
\end{defn0}

Note that $ i_{Z} $ is the maximum $ i $ such that $  Dh_{Z}(i+1)>0 $, i.e., the last element of the $h$-vector of $Z$. 
Additionally, if $ Z $ does not satisfy $\mathit{CB}(i)$, then it does not satisfy $\mathit{CB}(j)$ for 
all $ j \geq i $ either, by the contrapositive of Remark \ref{rem:CBprop}. Therefore, if $ Z $ satisfies $\mathit{CB}(i)$, then $ i \leq i_{Z} $.

\begin{exa0} Some examples of the Cayley-Bacharach property are shown below.
\begin{enumerate}[(i)]
 \item Let $ Z$ be a set of $6$ general points in $\Pj^{2} $. Then $ Dh_{Z} = (1,2,3) $, $ i_{Z} = 1 $ and $ Z $ satisfies $\mathit{CB}(1)$.
 \item Let $ Z$  be a set of $6$ general points lying on an irreducible conic of  $\Pj^{2} $.  Then $ Dh_{Z} = (1,2,2,1) $, $ i_{Z} = 2 $ and $ Z $ satisfies $\mathit{CB}(2)$, and, hence, $\mathit{CB}(1)$.
 \item Let $ Z$ be a set of $6$ general points in $\Pj^{2} $, with $ 5 $ of them on a line plus one point off the line. Then $ Dh_{Z} = (1,2,1,1,1) $, $ i_{Z} = 3 $ and $ Z $ does not satisfy $\mathit{CB}(1)$.
\end{enumerate}
 
\end{exa0}

The following holds.

\begin{lemma0}\label{lemma:tool}
Let $ Z = \{P_{1}, \ldots, P_{r}\} \subset \Pj^{n} $ be a finite set satisfying $\mathit{CB}(i)$. If, for any $ j \in \{1, \ldots, r\} $, 
the set $ Z_{j} = Z \setminus \{P_{j}\} $ does not satisfy $\mathit{CB}(i)$, then  
\begin{equation}\label{eq:h1Z}
h_{Z}(i+1) = \ell(Z)=r, \qquad\mbox{and so}\qquad Dh_Z(i+2)=0.
\end{equation}%
\end{lemma0}%
\begin{proof}
By hypothesis, for every $ j \in \{1, \ldots, r\} $ there exists a point $ Q_{j} \in Z_{j} $ and a form $ F_j \in H^{0}(\Oc_{\Pj^{n}}(i)) $ 
which vanishes at $ Z \setminus \{P_{j}, Q_{j}\} $ but not at $ Q_{j} $. Notice that $ F_j(P_{j})$ cannot be $0 $, 
since $ Z $ satisfies $\mathit{CB}(i)$. Let $ H_j \in H^{0}(\Oc_{\Pj^{n}}(1)) $ be a linear form vanishing at $ Q_{j} $ 
but not at $ P_{j} $ and consider $ G_j = F_j \cdot H_j $. For any $ j \in \{1, \ldots, r\} $, 
it holds that $ G_j \in H^{0}(J_{Z_{j}}(i+1)) \setminus H^{0}(J_{Z}(i+1)) $, which implies that $ Z $ does not satisfy 
$\mathit{CB}(i+1)$. Moreover, for all $j$, $ G_j(P_{j}) \neq 0  $ and $ G_j(P_{k}) = 0  $ for any $ k \neq j $. Therefore, in the exact sequence
$$ 0 \to H^{0}(J_{Z}(i+1)) \longrightarrow H^{0}(\Oc_{\Pj^{n}}(i+1)) \buildrel \rm {\it ev_{Z}(i+1)} 
\over \longrightarrow H^{0}(\Oc_{Z}(i+1))  $$
we have that $ \imm(ev_Z(i+1)) = H^{0}(\Oc_{Z}(i+1)) $, i.e., $ ev_{Z}(i+1) $ is a surjective map, which implies \eqref{eq:h1Z}.
\end{proof}

The following result, due to Geramita, Kreuzer, and Robbiano, gives a strong bound on the Hilbert function of sets with a Cayley-Bacharach property.

\begin{thm0}[Geramita, Kreuzer, and Robbiano \cite{GerKreuzerRobbiano93}]\label{CBprop}
If a finite set $ Z \subset \Pj^{n} $ satisfies $\mathit{CB}(i_{Z})$, then we have
\begin{enumerate}[(i)]
\item $ Dh_{Z}(0)+Dh_{Z}(1)+\cdots + Dh_{Z}(j) \leq Dh_{Z}(i_{Z}+1-j)+\cdots +Dh_{Z}(i_{Z}+1) $, 
for any $ j $ with $ 0 \leq j \leq i_{Z}+1 $;
\item $ Dh_{Z}(0)+Dh_{Z}(1)+\cdots + Dh_{Z}(j) \leq Dh_{Z}(k-j)+\cdots +Dh_{Z}(k) $, for any $ j, k $ 
with $ 0 \leq j \leq k \leq i_{Z}+1 $.
\end{enumerate}
\end{thm0}

\begin{proof} See Corollary 3.7 part (b) and (c) of \cite{GerKreuzerRobbiano93}.
\end{proof}

In order to create a link between Cayley-Bacharach properties and identifiability of symmetric tensors, 
we need to extend Theorem \ref{CBprop} by replacing $i_Z$ with any integer $i$ such that $Z$ satisfies $\mathit{CB}(i)$.

\begin{thm0}\label{GKRext}
If a finite set $ Z \subset \Pj^{n} $ satisfies $\mathit{CB}(i)$, then for any $ j $ such that $ 0 \leq j \leq i+1 $ we have
$$ Dh_{Z}(0)+Dh_{Z}(1)+\cdots + Dh_{Z}(j) \leq Dh_{Z}(i+1-j)+\cdots +Dh_{Z}(i+1).$$
\end{thm0}

\begin{proof} 
We proceed by induction on the residual part $h^{1}_{Z}(i)$, which is defined as:
$$ 
h^{1}_{Z}(i) 
= \ell(Z)-h_Z(i) 
= \ell(Z) - \sum_{j = 0}^{i} Dh_{Z}(j) 
= \sum_{j=i+1}^{i_{Z}+1} Dh_{Z}(j).
$$

If $ h^{1}_{Z}(i) = 1 $, then $\ell(Z)=h_Z(i)+1$. Thus, by Remark \ref{resto1} we must have $\ell(Z)=h_Z(i+1)$, i.e., $ i = i_{Z} $  and we conclude by Theorem \ref{CBprop} part (i).

Next, we assume that the theorem is true for all $1 \le h^{1}_{Z}(i) \le e$, and prove it for $e+1$.
If $ Z $ satisfies $\mathit{CB}(i_{Z})$, then we can conclude by Theorem \ref{CBprop} part (ii). 
So assume that $ Z $ does not satisfy $\mathit{CB}(i_{Z})$. In this case notice that
there exists $ P \in Z $ such that $ Z_{P} = Z \setminus \{P\} $ satisfies $\mathit{CB}(i)$. Indeed,
if for any $ P \in Z $ the set $ Z_{P} $ does not satisfy $\mathit{CB}(i)$, then by Lemma \ref{lemma:tool} 
we have that $ h^{1}_{Z}(i+1) = 0 $, that is $ i = i_{Z} $, which is a contradiction. 
Fix then a point $ P \in Z $ such that $ Z_{P} $ satisfies $\mathit{CB}(i)$. 
Since $ Z $ satisfies  $\mathit{CB}(i)$, by \eqref{eq:h0} we have that 
$$ h_{Z_{P}}(j) = h_{Z}(j)\qquad  \mbox{ and }\qquad Dh_{Z_{P}}(j) = Dh_{Z}(j) $$
for $ j \leq i $. 
Therefore,
$$ h^{1}_{Z_P}(i)=\ell(Z_P)- \sum_{j = 0}^{i} Dh_{Z_{P}}(j) = \ell(Z)-1 -\sum_{j = 0}^{i} Dh_{Z}(j) =h^{1}_{Z}(i)-1, $$
so that by applying the induction hypothesis to $ Z_{P} $, we get
$$ Dh_{Z_{P}}(0)+Dh_{Z_{P}}(1)+\cdots + Dh_{Z_{P}}(j) \leq Dh_{Z_{P}}(i+1-j)+\cdots +Dh_{Z_{P}}(i+1) $$
for every $ j $ such that $ 0 \leq j \leq i+1 $. From \eqref{eq:h0} and $ Dh_{Z_{P}}(i+1) \leq Dh_{Z}(i+1) $, by \reflem{rem:triv} (viii), we get the conclusion.
\end{proof}

\section{Castelnuovo's Lemma revisited}\label{sec:Crev}

In order to apply the Cayley-Bacharach property to the identifiability of symmetric tensors, we need an {\it ad hoc} 
extension of the classical Castelnuovo Lemma for finite sets in $\Pj^n$.

\begin{lemma0}[Castelnuovo, \cite{GH} page 531]
Let $ Z \subset \Pj^n $ be a finite set such that:
\begin{enumerate}[(i)]
\item $ \ell(Z) \geq 2n+3 $;
\item $ Z $ is in LGP;
\item $ Z $ imposes $ 2n+1 $ conditions to quadrics, i.e., $ h_{Z}(2) = 2n+1 $. 
\end{enumerate}
Then, the quadrics containing $ Z $ intersect in a rational normal curve, which thus contains $Z$.
\end{lemma0}

In particular, we need to weaken the hypothesis, by assuming only that some
subset of $Z$ is in LGP. We do that in several steps.

\begin{prop0}\label{prop:intermediatestep}
Let $ Z \subset \Pj^n $ be a finite set with a subset $ A \subset Z $ of cardinality $ 2n+1 $ in LGP. Then $ h_{Z}(2) \geq 2n+1 $.
\end{prop0}
\begin{proof}
Since $ A $ is in LGP, then $ A $ imposes independent conditions to quadrics. Therefore $ h_{A}(2) = 2n+1 $. As $ A \subset Z $, \reflem{rem:triv} (viii) concludes the proof.
\end{proof}

By means of Proposition \ref{prop:intermediatestep}, Castelnuovo's Lemma can be rephrased as follows:

\begin{lemma0}\label{lemma:Cclass}
 Let $ Z \subset \Pj^n $ be a finite set such that:
\begin{enumerate}[(i)]
\item $ \ell(Z) \geq 2n+3 $;
\item $ Z $ is in LGP;
\item $ Z $ imposes at most $ 2n+1 $ conditions to quadrics, i.e., $ h_{Z}(2) \leq 2n+1 $. 
\end{enumerate}
 Then, $ h_{Z}(2) = 2n+1 $ and the quadrics containing $ Z $ intersect in a rational normal curve, which thus contains $Z$.
\end{lemma0}

Based on Lemma \ref{lemma:Cclass}, we can prove the following extension:

\begin{lemma0}\label{lemma:Cnew}
 Let $ Z \subset \Pj^n $ be a finite set such that:
\begin{enumerate}[(i)]
\item $ \ell(Z) \geq 2n+3 $;
\item $ h_{Z}(2) \leq 2n+1 $;
\item there exists $ A \subset Z $ such that $ \ell(A) = 2n+1 $ and $ A $ is in LGP. 
\end{enumerate}
Then, $ Z $ is in LGP and it is contained in a rational normal curve.
\end{lemma0}
\begin{proof} First, let us assume that $ \ell(Z) = 2n+3 $ and let us set $ Z = A \cup \{P,Q\} $. We 
claim that $ W = A \cup \{P\} $ is in LGP. 
Indeed, if this is not the case, there are subsets of $W$, of cardinality at most $n+1$, 
which are not linearly independent. This implies that one can find $ V \subset W $ and a hyperplane $ H $ such that 
$ \ell(V) = n+1 $ and $ V \subset H $. Since $ A $ is in LGP, then, necessarily, $ V $ contains $ P $ and 
 $ n $ points of $ A $. Thus we can renumber the points
of $ A = \{P_1,\dots,P_{2n+1}\} $ so that $ P_i\in H $ if and only if $i \leq n $.
For any $ j \in \{n+1,\dots, 2n \} $, we can find a hyperplane $ H_{j} \subset \Pj^{n} $ such that 
$ \{P_{n+1},\dots, P_j\} \subset H_{j} $ and $ \{P_{j+1},\dots, P_{2n+1}\} \not\subset H_{j} $, 
because $ A $ is in LGP. It follows that the quadric $ Q_j=H\cup H_j $
contains $ \{P,P_1,\dots, P_j \}$ and misses $ \{P_{j+1},\dots, P_{2n+1}\} $.
In particular, if we set $ V_{j} = \{P,P_{1}, \ldots,P_{j}\} $, then $ h^{0}(J_{V_{j}}(2)) > h^{0}(J_{V_{j+1}}(2)) $,
which implies:
\begin{equation}\label{eq:h2}
h_{V_{j}}(2) <  h_{V_{j+1}}(2).  
\end{equation}
We notice that $ h_{V_{j}}(1) = n+1 $, for $V_j$ contains at least $n+1$ points of $ A $, which
 is in LGP. Therefore $ Dh_{V_{j}}(1) = n $. Moreover, since $\ell(V_{n+1})=n+2=h_{V_{n+1}}(1)+1$, 
 by Remark \ref{resto1} we get that $h_{V_{n+1}}(2)=\ell(V_{n+1})$. Thus, the $h$-vector of $V_{n+1}$
 is $(1,n,1)$.
By induction on $ j $, we show that
\begin{equation}\label{eq:ind}
h_{V_{j}}(2) = j+1,\quad \forall j\geq n+1. 
\end{equation}
Indeed, the claim holds for $ j = n+1 $. If \eqref{eq:ind} holds for $ j \geq n+1 $, 
then by  \eqref{eq:h2} we have that $ h_{V_{j+1}}(2) > j+1 $. Since, by definition, $ h_{V_{j+1}}(2) \leq \ell(V_{j+1}) = j+2 $, necessarily it has to be the case that $ h_{V_{j+1}}(2) = j+2 $, as desired.
By applying \eqref{eq:ind} to $ j = 2n+1 $, we get that 
$$ h_{W}(2) = h_{V_{2n+1}}(2) =  2n+2, $$
which contradicts assumption (ii) via \reflem{rem:triv} (viii), as $ W \subset Z $. So, $ W $ is in LGP.

Now assume that $Z$ is not in LGP. Then, as above, there exists $ V \subset Z $ and a hyperplane 
$ H \subset \Pj^{n} $ such that $ \ell(V) = n+1 $ and $ Q \in V \subset H $. Then there exists $ U \subset Z $ 
such that $ \ell(U) = 2n+2 $ and $ U $ is not in LGP. Since $ h_{U}(2) \leq 2n+1 = h_Z(2)$, and $ U $ contains $ 2n+1 $ 
points of $ W $ which is in LGP, then we get a contradiction by arguing as above.

Finally, assume that $\ell(Z)>2n+3$. Notice that we have just proved the existence of $ Z_{0} \subset Z $ 
such that $ \ell(Z_{0}) = 2n+3 $ and $ Z_{0} $ is in LGP. In particular, $ h_{Z_{0}}(2) \leq 2n+1 $. 
By Lemma \ref{lemma:Cclass}, $ h_{Z_{0}}(2) = 2n+1 $ and $ Z_{0} $ lies in a rational normal curve $\Gamma$, 
which is the intersection of all quadrics containing $Z_{0}$. Since $h_Z(2)\leq 2n+1$, then equality holds and 
$ H^{0}(J_{Z_{0}}(2)) = H^{0}(J_{Z}(2)) $. Thus $ Z $ itself is contained in $\Gamma$, which, 
by Remark \ref{rem:rncLGP},  implies that $ Z $ is in LGP.
\end{proof}

\begin{exa0}
The previous formulation of Castelnuovo's Lemma is sharp, in the sense that the existence of
a subset of cardinality $2n$ in LGP is not enough to guarantee that a set $Z$ of $2n+3$ points in $\Pj^n$,
with $h_Z(2)=2n+1$, is contained in a rational normal curve.

Namely, take $n=3$ and take a smooth quadric $Q \subset \Pj^3$, a set of $4$ general points $Q_1,\dots,Q_4$ on $Q$, a general line
$L\subset Q$ of type $(0,1)$ and a set of $5$ general points $P_1,\dots,P_5$ on $L$. The set $Z=\{Q_1,\dots,Q_4,P_1,\dots,P_5\}$
 has cardinality $9=2n+3$ and contains the subset $A=\{Q_1,\dots,Q_4,P_1,P_2\}$ of cardinality
 $6=2n$ which is in LGP,  due to the generality in the choice of the points.
 Since $4$ general points of $Q$ lie in two linearly independent divisors of type $(2,1)$, then there are $3$
 independent quadrics of $\Pj^3$ containing $Z$, i.e., $h_Z(2)=7=2n+1$. The set $Z$, however, cannot lie
 in a rational normal curve, for it contains $5$ points on a line.
 \end{exa0}

\section{Application to the identifiability of quartics}\label{sec_application}

In this section, we apply the previous results on the geometry of finite sets to decompositions
of quartic polynomials with the purpose of reaching beyond the range of Kruskal's criterion of identifiability.

Throughout this section, we consider $d=4$ and let $ T \in S^{4}\mathbb{C}^{n+1} $ be a homogeneous polynomial.
Let $ A \subset \Pj^n $ be a finite set that computes $T$, such that:
\begin{enumerate}[(i)]
\item $ \ell(A) = 2n+1 $, i.e., $A$ has one point more than Theorem \ref{thm:Kruskalresh} allows; see Remark \ref{rem:last};
\item $ A $ is in LGP;
\item $ A $ satisfies the minimality property.
\end{enumerate}
Let $ B \subset \Pj^n $ be \textit{another} finite set that computes $T$ with $ \ell(B) \leq 2n+1 $. Without loss of generality we can also assume that $ B $ satisfies the minimality property.
In particular, $ T \in \langle \nu_{4}(A) \rangle \cap \langle \nu_{4}(B) \rangle $.
Our target is to find criteria that exclude the existence of $B$, so that $T$ has rank $2n+1$
and is identifiable. 

In the following, we analyze the geometry of the union
$$ Z=A\cup B.$$
Clearly, $\ell(Z)\leq \ell(A)+\ell(B)$, with equality if $A\cap B$ is empty.
%
It is a straightforward fact that
$$\langle \nu_{4}(Z)\rangle=\langle \nu_{4}(A)\rangle+ \langle \nu_{4}(B)\rangle.$$
so that by using Grassmann's formula, we have that
\begin{equation}\label{eq:vZ1}
\dim(\langle \nu_{4}(Z)\rangle)=\dim(\langle \nu_{4}(A)\rangle)+
\dim(\langle \nu_{4}(B)\rangle)-\dim(\langle \nu_{4}(A)\rangle\cap\langle \nu_{4}(B)\rangle).
\end{equation}
%
%
From Remark \ref{rem:indep}, we then find that
\begin{equation}\label{eq:Gr}
\dim(\langle \nu_{4}(Z)\rangle) = \ell(A)+\ell(B)-2- \dim(\langle \nu_{4}(A)\rangle\cap\langle \nu_{4}(B)\rangle).
\end{equation}
Since $A\cap B $ is a proper subset of $ A$, then, by the minimality assumption on $ A $, we get that 
$ T \not\in \langle \nu_{4}(A\cap B)\rangle$. Thus,
\begin{equation}\label{eq:int}
\dim(\langle \nu_{4}(A)\rangle\cap \langle \nu_{4}(B)\rangle) > \dim(\langle \nu_{4}(A\cap B)\rangle).  
\end{equation}
From (\ref{eq:Gr}), (\ref{eq:int}) and since $ \dim(\langle \nu_{4}(A\cap B)\rangle) = \ell(A\cap B)-1 $, 
it follows that
\begin{align}\label{eq:h1e}
\dim(\langle \nu_{4}(Z)\rangle) &< \ell(A)+\ell(B)-2- \dim(\langle \nu_{4}(A\cap B)\rangle) \\
\nonumber &= \ell(A)+\ell(B)-\ell(A\cap B)-1 =\ell(Z)-1. 
\end{align}
%
We notice that
\begin{equation*}\label{eq:h4}
h_Z(4)\leq \ell(Z)-1.
\end{equation*}
Indeed by the inequality (\ref{eq:h1e}), the dimension
of $\langle \nu_{4}(Z)\rangle$, which is, by Remark \ref{rem:conds}, $h_Z(4)-1$, cannot be $\ell(Z)-1$. Thus, by \refprop{nonincr} and \reflem{rem:triv} (v):
\begin{equation}\label{eq:Dh5}
Dh_Z(5)>0.
\end{equation}
%
We set now, as usual,  
$$h^1_Z(4)=\ell(Z)-h_Z(4)=\sum_{j=5}^\infty Dh_Z(4).$$
%
By means of (\ref{eq:h^1}) and part (iv) of \reflem{rem:triv}, we have that
\begin{align}\label{eq:vZ}
\dim(\langle \nu_{4}(Z)\rangle)
&= h_Z(4)-1 =\ell(Z)-1-h^1_Z(4) \\
\nonumber &= \ell(A)+\ell(B)-\ell(A\cap B)-1-h^1_Z(4).
\end{align}
Since $\ell(A)=\dim(\langle \nu_{4}(A)\rangle)+1$ and $\ell(B)=\dim(\langle \nu_{4}(B)\rangle)+1$, 
then, comparing (\ref{eq:vZ}) and (\ref{eq:vZ1}) we get that
\begin{equation}\label{h1span}
\dim(\langle \nu_{4}(A)\rangle\cap\langle \nu_{4}(B)\rangle) =\ell(A\cap B)-1+h^1_Z(4).
\end{equation}

\begin{prop0}\label{prop:lenghtZ}
We cannot have that $\ell(Z) \leq \ell(A)+1=2n+2$.
\end{prop0}
\begin{proof} From Proposition \ref{prop:intermediatestep} and our assumptions, we know that  
$\ell(A)=h_A(2)$, i.e., quadrics separate the points of $A$. Thus $ h_{Z}(2) \geq  h_A(2) = 2n+1 $. 
It follows that if $\ell(Z) \leq 2n+2$, then  $ h_{Z}(2) \geq \ell(Z)-1  $, 
which implies, by Remark \ref{resto1}, that $ h_Z(3)=\ell(Z)$; thus, $Dh_Z(4) = 0 $. This fact contradicts
Proposition \ref{nonincr}, since (\ref{eq:Dh5}) holds. 
\end{proof}

We use the previous arguments and Lemma \ref{lemma:Cnew} to prove the next crucial result.

\begin{thm0}\label{thm:main} 
With our assumptions, $ Z $ is contained in a rational normal curve.
\end{thm0}
\begin{proof} 
First, note that $ Z $ satisfies Lemma \ref{lemma:Cnew} (iii) and, because of \refprop{prop:lenghtZ} it also satisfies \reflem{lemma:Cnew} (i).
Assume that $Z$ has the Cayley-Bacharach property $\mathit{CB}(4)$. By applying Theorem \ref{GKRext} 
with $ i = 4 $ and $ j = 2 $ we get
\begin{equation}\label{eq:CB4}
Dh_Z(0) + Dh_Z(1) + Dh_Z(2) \leq Dh_Z(3)+Dh_Z(4)+Dh_Z(5).
\end{equation}
Since $Z$ contains $n+1$ points of $A$ in LGP, then $ \langle Z\rangle = \Pj^n $, so that $Dh_Z(1)=n$. 
We claim that also (ii) of Lemma \ref{lemma:Cnew} holds. Indeed, if this is not the case, then $Dh_Z(2)>n $
and thus, by \eqref{eq:CB4} and \reflem{rem:triv} part (v),
$$ \ell(Z) = \sum_j Dh_Z(j) \geq \sum_{j=0}^{5} Dh_Z(j) > 4n+2.  $$
On the other hand, $\ell(Z)\leq \ell(A) + \ell(B) \leq 2\ell(A) = 4n+2 $, which leads to a contradiction. 
Therefore, all the assumptions of Lemma \ref{lemma:Cnew} hold, concluding the proof for this case.

It remains to prove that with our assumptions $Z$ necessarily has the property $\mathit{CB}(4)$.
 Assume that $Z$ does not satisfy $\mathit{CB}(4)$. Then, there exists $P\in Z$
such that $ h^{0}(J_{Z\setminus\{P\}}(4)) > h^{0}(J_{Z}(4)) $, i.e., $ h_{Z\setminus\{P\}}(4) < h_Z(4) $. 
As $ \ell(Z\setminus\{P\}) = \ell(Z) - 1 $, by \reflem{rem:triv} (v), we have $ Dh_{Z\setminus\{P\}}(j) 
< Dh_Z(j) $ for some $ j \in \{1,2,3,4\} $. Since $ Dh_{Z\setminus\{P\}}(j) \leq Dh_Z(j) $ for all $ j $, it follows that 
$ Dh_{Z\setminus\{P\}} (j) = Dh_Z(j) $
for $ j\geq 5 $, i.e., by \eqref{eq:h^1}, $ h^1_Z(4) = h^1_{Z\setminus\{P\}}(4)$.

Assume that $ A \cap B = \emptyset $. If $ P\in A $, then, by \eqref{h1span}, $\dim (\langle \nu_{4}(A\setminus\{P\})\rangle
\cap\langle \nu_{4}(B) \rangle) = \dim (\langle \nu_{4}(A)\rangle\cap\langle \nu_{4}(B)\rangle) $,
so $ \langle \nu_{4}(A\setminus\{P\})\rangle\cap\langle \nu_{4}(B) \rangle = \langle \nu_{4}(A)\rangle\cap
\langle \nu_{4}(B)\rangle $. Thus $ A\setminus\{P\} $ computes $ T $, which contradicts the hypothesis of minimality of $ A $. 
Therefore, $ P\in B\setminus A $. But now we can repeat the previous argument for $ A $ and $ B \setminus\{P\} $,
and we get that $B$ is not minimal, which is a contradiction. 

Hence, we can assume that $ A \cap B \not = \emptyset $. We claim that
\begin{equation}\label{eq:int2}
\langle \nu_{4}(A) \rangle \cap \langle \nu_{4}(B \setminus A) \rangle \not = \emptyset.  
\end{equation}
Indeed, let $ s = \ell(A \cap B) $. We can renumber the elements of $  \nu_{4}(A) $ and $  \nu_{4}(B) $ 
in a way such that, for both sets, the first $ s $ comprise $ \nu_{4}(A \cap B) $, i.e.,
$$ \nu_{4}(A \cap B) = \{\nu_{4}(P_{1}), \ldots, \nu_{4}(P_{s})\}. $$ 
Note that $ \langle \nu_{4}(A \cap B) \rangle $ is a proper subset of $ \langle \nu_{4}(A) \rangle \cap 
\langle \nu_{4}(B) \rangle $, since, for example, $ T \in (\langle \nu_{4}(A) \rangle \cap \langle \nu_{4}(B) \rangle)
 \setminus \langle \nu_{4}(A \cap B) \rangle $, as $ A $ and $ B $ are minimal for $ T $. Therefore, 
\begin{align*}
T 
&= \alpha_{1} \nu_{4}(P_{1})+ \ldots + \alpha_{s} \nu_{4}(P_{s}) + \alpha_{s+1} \nu_{4}(P_{s+1}) + 
\ldots + \alpha_{\ell(A)} \nu_{4}(P_{\ell(A)}) \\
&= \beta_{1} \nu_{4}(P_{1})+ \ldots + \beta_{s} \nu_{4}(P_{s}) + \beta_{s+1} \nu_{4}(Q_{s+1}) + 
\ldots + \beta_{\ell(B)} \nu_{4}(Q_{\ell(B)}) 
\end{align*}
with $ \alpha_{j}, \beta_{k} \in \C \setminus\{0\} $ for any $ j, k $ and $ P_{j} \in A \setminus B, 
\, Q_{k} \in B \setminus A $ for $ j \in \{s+1, \ldots, \ell(A)\} $ and $ k \in \{s+1, \ldots, \ell(B)\} $. It turns out that
$$ \sum_{i=1}^{s}(\alpha_{i} - \beta_{i}) \nu_{4}(P_{i}) + \sum_{i=s+1}^{\ell(A)}\alpha_{i} \nu_{4}(P_{i}) - 
\sum_{i=s+1}^{\ell(B)} \beta_{i} \nu_{4}(Q_{i}) ) = 0. $$
Thus, the tensor
\begin{equation}\label{eq:twodec}
T' = \sum_{i=1}^{s}(\alpha_{i} - \beta_{i}) \nu_{4}(P_{i}) + \sum_{i=s+1}^{\ell(A)}\alpha_{i} \nu_{4}(P_{i}) 
= \sum_{i=s+1}^{\ell(B)} \beta_{i} \nu_{4}(Q_{i}) )
\end{equation}
is an element of $ \langle \nu_{4}(A) \rangle \cap \langle \nu_{4}(B \setminus A) \rangle $, which implies \eqref{eq:int2}. 
Now, if $ A $ is minimal for $ T' $, then we have two finite sets $ A' = A $ and $ B' = B \setminus A $ computing $ T' $ 
and such that $ A' \cap B' = \emptyset $. Thus, by replacing $ T $ with $ T' $ and by arguing as in the case 
$ A \cap B = \emptyset $, we get a contradiction because $A'\cup B'=A\cup B$ does not satisfy $\mathit{CB}(4)$. 
Therefore, we can assume that $ A  $ is not minimal for $ T' $, i.e., there exists a proper subset $ A' $ of $ A $ such that 
$ T' \in \langle\nu_{4}(A') \rangle $. Then some of the $ P_{i} $'s with $ i \in \{1, \ldots, s\} $ does not appear in the decomposition 
of $ T' $, say $ P_{1} $. So there exists $ \gamma_{i} \in \C \setminus \{0\} $,  such that
$$ T' = \sum_{i=1}^{s}(\alpha_{i} - \beta_{i}) \nu_{4}(P_{i}) + \sum_{i=s+1}^{\ell(A)}\alpha_{i} \nu_{4}(P_{i}) = 
\gamma_{2}\nu_{4}(P_{2}) +\ldots +\gamma_{\ell(A)} \nu_{4}(P_{\ell(A)}). $$
Since, by Remark \ref{rem:indep}, $ \nu_{4}(P_{1}), \ldots, \nu_{4}(P_{\ell(A)}) $ are linearly independent, 
it follows that $ \alpha_{1} = \beta_{1} $, $ \gamma_{i} = \alpha_{i} - \beta_{i} $ for $ i \in \{2, \ldots, s\} $ 
and $ \gamma_{i} = \alpha_{i} $ for $ i \in \{s+1, \ldots, \ell(A)\} $. Therefore, by \eqref{eq:twodec}, 
$$ T' = \sum_{i=2}^{\ell(A)} \gamma_{i} \nu_{4}(P_{i}) = \sum_{i=s+1}^{\ell(B)} \beta_{i} \nu_{4}(Q_{i}) ), $$  
so that $ T' $ has two different decompositions $A',B'$ with, respectively, $\ell(A') = \ell(A) - 1 = 2n $ and $  \ell(B') = \ell(B) - s \leq 2n $ summands. 
As $\ell(A'),\ell(B')\leq 2n$, this contradicts Theorem \ref{thm:Kruskalresh}.

\end{proof}

As a consequence of Theorem \ref{thm:main} we get the following result.

\begin{thm0}\label{main2} Fix a homogeneous polynomial $ T \in S^{4}\mathbb{C}^{n+1} $ for which there exists a finite,
minimal set $A$ that computes $T$, such that $\ell(A)=2n+1$ and $A$ is in LGP.
Then, the existence of a second set $B\subset\Pj^n$ that computes $T$ and $\ell(B)\leq 2n+1$ implies that $\ell(B)=2n+1$
and both $A$, $B$ belong to a rational normal curve of $\Pj^n$. 
\end{thm0}

If we want to understand the identifiability of quadrics of rank $2n+1$, we should thus study the case where a minimal set that computes $T$ lies in a rational normal curve.

\begin{prop0} \label{prop_criterion}
Fix a homogeneous polynomial $ T \in S^{4}\mathbb{C}^{n+1} $ for which there exists a finite,
minimal set $A$ that computes $T$, such that $\ell(A)=2n+1$ and $A$ is contained in a rational normal curve 
$\Gamma\subset\Pj^n$.
Then, there exists a positive dimensional family of finite sets $ A_{t} \subset \Gamma $ such that 
\begin{enumerate}[(i)]
\item $ \ell(A_{t}) = 2n+1 $;
\item $ A_{0} = A $;
\item $ T \in \langle \nu_{4}(A_t) \rangle$, i.e., each $A_t$ computes $T$.
\end{enumerate}
\end{prop0}
\begin{proof} The curve $\Gamma$ is the image of a Veronese map $ \Gamma = \nu_{n}(\Pj^1) $, thus 
$ \Gamma' = \nu_{4}(\Gamma) = \nu_{4n}(\Pj^1) $ is a rational 
normal curve in $ \Pj^{4n} $. Therefore, $ T $ belongs to the secant variety $ \Sec{2n+1}{\Gamma'} $. 
By Remark \ref{rem:rncsec}, $ \Sec{2n+1}{\Gamma'} $ covers $ \Pj^{4n} $. 
If $ \Sigma_{2n+1}(\Gamma') $ denotes the abstract $(2n+1)$-secant variety of $ \Gamma' $, 
then the fibre of the $ (2n+1)$-secant map
$$ \pi_{2n+1} : \Sigma_{2n+1}(\Gamma') \rightarrow \Pj^{4n} $$
 at $ T $ has dimension $ \dim(\Sigma_{2n+1}(\Gamma')) - 4n = 2(2n+1) - 1- 4n = 1 $. 
\end{proof}

\begin{rem0}\label{rem:dep1}
With the above notation, the tangent lines to $ \Gamma' $ at the points of $ \nu_{4}(A) $ span a space of 
dimension at most $ 4n $. Therefore, the tangent spaces to $\nu_{4}(\Pj^n)$ at the same points span a space 
of dimension at most $ (2n+1)(n+1)-2 $.
\end{rem0}

Summarizing, for quartics $T$ in $\Pj^n$ which are computed by a set $A\subset\Pj^n$ in LGP and
cardinality at most $2n+1$, we have that either:
\begin{enumerate}[(i)]
\item $T$ has rank $\ell(A)$ and is identifiable; or
\item $T$ has rank $2n+1$ and is computed by a $1$-dimensional family of sets of cardinality $2n+1$ that includes $A$.
\end{enumerate}

One can use \textit{Terracini's test}, introduced in \cite[Lemma 6.5]{COttVan17b}, to decide which of these two cases occurs.


\subsection{The algorithm} \label{sec_the_test}
In an abuse of notation, let $\nu_4 : \C^{n+1} \to \C^{\binom{n+4}{4}}$ be the $4$-fold symmetric tensor product. Then, $[\nu_4(\vect{m})] = \nu_4([\vect{m}])$.
Given a length-$r$ symmetric tensor rank decomposition of a quartic 
\[
 T = \sum_{i=1}^r \nu_{4}( P_i )
\]
in the form of the collection of points $A = \{ P_i = [\vect{m}_i] \}_{i=1}^r \subset \Pj^n$, we can apply the following algorithm for verifying that the given decomposition of $T$ is identifiable:
\begin{enumerate}
 \item[S1.] If $r > 2n + 1$, the criterion cannot be applied.
 \item[S2.] If $r < 2n + 1$, use the reshaped Kruskal criterion from \cite[Section 6.2]{COttVan17a}.
 \item[S3.] If $r = 2n + 1$, perform the next tests: 
 \begin{enumerate}
  \item[1)] \emph{minimality test}: check that $\dim \langle \nu_4(\vect{m}_1), \ldots, \nu_4(\vect{m}_r) \rangle = r$;
  \item[2)] \emph{Kruskal's test}: check that $A$ is in LGP;
  \item[3)] \emph{Terracini's test}: check that $\dim \langle \Tang{\vect{m}_1}{\nu_4(\C^{n+1})}, \ldots, \Tang{\vect{m}_r}{\nu_4(\C^{n+1})} \rangle = 2n^2 + 3n + 1$.
 \end{enumerate}
 If all these tests are successful, then $T$ is of rank $r$ and is $r$-identifiable.
\end{enumerate}

An implementation of this algorithm is included in the ancillary Macaulay2 file \texttt{identifiabilityS4Cn.m2}.


We note that the new criterion for $r = 2n+1$ is effective in the sense of \cite{COttVan17b}. Indeed, quartics with $r = 2n+1$ are always generically $r$-identifiable \cite{Ball05a}, and it is easy to verify that the conditions in tests $1$, $2$, and $3$ are not satisfied precisely on a Zariski-closed strict subvariety of the $r$-secant variety of $\nu_4(\Pj^n)$.

\subsection{Examples}
We present some examples of identifiable and unidentifiable Waring decompositions in the original case $r = 2n+1$.

\paragraph{\textit{An identifiable example}} 
Consider $n = 4$ and $r = 2n+1 = 9$. We generated a random collection of $9$ points $A = \{P_i = [\vect{m}_i] \}_{i=1}^9$ in Macaulay2, where the vectors $\vect{m}_i \in \N^5$ had the following values in our experiment
\[
M = \begin{bmatrix}
 \vect{m}_i
\end{bmatrix}_{i=1}^9 = 
\begin{bmatrix}
 0  &  1 &  1 & -3 & -5 &  2 & -1 &  2 & -1 \\
 -2 & -1 &  2 &  0 &  1 &  2 & -4 &  3 & 1  \\
 2  &  0 &  5 &  1 &  4 & -5 & -1 & -3 & 4  \\
 1  & -5 & -1 &  3 & -2 &  3 &  5 &  2 & -3 \\
 1  & -3 & -2 & -5 & -4 &  3 & -2 &  1 & 4  
\end{bmatrix}.
\]
The minimality test shows that $\dim \langle \nu_4(A) \rangle = \operatorname{rank}([\vect{m}_i \otimes \vect{m}_i \otimes \vect{m}_i \otimes \vect{m}_i]_{i=1}^9) = 9$, which is as required. We then compute the rank of all $126$ subsets of $5$ columns of $M$. They are all of rank $5$, so that the Kruskal rank is $5$ and $A$ is in LGP. Finally, we compute a basis $B_i$ of the affine cone over the tangent space of the Veronese variety $\nu_4(\Pj^n)$ at one of the points in the cone over $\nu_4(P_i)$, and then compute $\dim \langle B_1, \ldots, B_9 \rangle$. The computation reveals that it equals $45 = 2 \cdot 4^2 + 3\cdot4 +1$, so that Terracini's test is also successful. We can conclude that $T = \sum_{i=1}^9 \nu_4(\vect{m}_i)$ is (complex) identifiable.

\paragraph{\textit{A variation of Derksen's example}} 
Consider the Vandermonde matrix 
\[
 M = [\vect{m}_i]_{i=1}^{r} = \begin{bmatrix}
      1 & 1 & \cdots & 1 \\
      \lambda_1 & \lambda_2 & \cdots & \lambda_r \\
      \lambda_1^2 & \lambda_2^2 & \cdots & \lambda_r^2 \\
      \vdots & \vdots & & \vdots \\
      \lambda_1^{n} & \lambda_2^n & \cdots & \lambda_r^n
     \end{bmatrix}.
\]
Let $A = \{ P_i = [\vect{m}_i]_i \}$ be the corresponding collection of points in $\Pj^n$. By construction, each of the points is in the image of $\nu_n( \Pj^1 )$, i.e., they lie on a rational normal curve in $\Pj^n$. Hence, if $\lambda_1,\ldots,\lambda_r \in \C$ are pairwise distinct then $A$ is in LGP. Taking $r = 2n+1$, then $T = \sum_{i=1}^r \nu_4(\vect{m}_i)$ is not symmetric identifiable by \refprop{prop_criterion}, even though $A$ is LGP. 

We applied the criterion to the case $n = 4$, $r = 9$, and $\lambda_i = i-1$ for $i=1,\ldots,9$. Running the algorithm, we find that both the minimality test and the Kruskal test are successful, consistent with the theory in \cite{Derksen13}. Our theory predicts that Terracini's test must fail, because $T$ can now only be unidentifiable if there is a rational normal curve passing through the points $A$. Performing the computation, we find that the tangent spaces only span a space of dimension $44$, one less than expected. Hence, Terracini's test fails, and we cannot conclude that $T$ is identifiable.

\bibliographystyle{amsplain}
\bibliography{biblioLuca}

\end{document}